\documentclass[a4paper,11pt,reqno]{amsart}

\usepackage{amssymb}

\newtheorem*{theorem*}{Theorem}

\theoremstyle{definition}

\newtheorem*{example*}{Example}

\newtheorem*{remark*}{Remark}

\newcommand{\bF}{{\mathbb F}}
\newcommand{\bZ}{{\mathbb Z}}

\newcommand{\calJ}{{\mathcal J}}

\newcommand{\calS}{{\mathcal S}}

\newcommand{\calL}{{\mathcal L}}
\newcommand{\calK}{{\mathcal K}}

\newcommand{\subo}{_{\bar 0}}
\newcommand{\subuno}{_{\bar 1}}

\providecommand{\espan}[1]{\text{span}\left\{ #1\right\}}

 \DeclareMathOperator{\frsl}{{\mathfrak{sl}}}
 
 \DeclareMathOperator{\frso}{{\mathfrak{so}}}

\begin{document}

\title{More non semigroup Lie gradings}

\author[Alberto Elduque]{Alberto Elduque$^{\star}$}
 \thanks{$^{\star}$ Supported by the Spanish Ministerio de
 Educaci\'{o}n y Ciencia
 and FEDER (MTM 2007-67884-C04-02) and by the
Diputaci\'on General de Arag\'on (Grupo de Investigaci\'on de
\'Algebra)}
\address{Departamento de Matem\'aticas e
Instituto Universitario de Matem\'aticas y Aplicaciones,
Universidad de Zaragoza, 50009 Zaragoza, Spain}
\email{elduque@unizar.es}

\date{\today}



\begin{abstract}
This note is devoted to the construction of two very easy examples, of respective dimensions $4$ and $6$, of graded Lie algebras whose grading is not given by a semigroup, the latter one being a semisimple algebra. It is shown that $4$ is the minimal possible dimension.
\end{abstract}

\maketitle



Patera and Zassenhaus \cite{PZ} define a \emph{Lie grading} as a
decomposition of a Lie algebra into a direct sum of subspaces
\[
\calL=\oplus_{g\in G}\calL_g,
\]
indexed by elements $g$ from a set $G$, such that $\calL_g\ne 0$ for
any $g\in G$, and for any $g,g'\in G$, either
$[\calL_g,\calL_{g'}]=0$ or there exists a $g''\in G$ such that
$0\ne [\calL_g,\calL_{g'}]\subseteq \calL_{g''}$.

Then, in \cite[Theorem 1.(d)]{PZ}, it is asserted that, given a Lie
grading, 
the set $G$ embeds in an abelian semigroup
so that the following property holds:
\begin{itemize}
\item[(P)]
For any $g,g',g''\in G$ with $0\ne [\calL_g,\calL_{g'}]\subseteq
\calL_{g''}$, $g+g'=g''$ holds in the semigroup.
\end{itemize}

In \cite{E}, a counterexample to this assertion was given. A nilpotent Lie algebra of dimension $16$ was defined, with a grading not given by an abelian semigroup. This example came as a surprise (see \cite{Svobodova}), but its difficulty may give the impression that this is a rare phenomenon.

The purpose of this note is to give two more counterexamples that show that non semigroup gradings are not so rare. The first one will be a non semigroup grading on a four dimensional solvable (actually metaabelian) Lie algebra. It will be shown that there are no counterexamples in dimension $\leq 3$, so this is a counterexample of minimal dimension. On the other hand, another such non semigroup grading will be defined over the direct sum of two three dimensional simple Lie algebras. It has the interesting feature of being a coarsening of a group grading. Note that in the last paragraph in \cite{E} a grading on $\frsl(2)\oplus\frsl(2)$ was shown which is not a group grading, however it is a semigroup grading.

\smallskip

Anyway, the question posed in \cite{E} still remains open:

\noindent \emph{Is any grading on a simple finite dimensional complex Lie algebra a group grading?}

\smallskip

Some notation is in order before we start.

Given a graded Lie algebra $\calL=\oplus_{g\in G}\calL_g$, a subspace $\calS$ of $\calL$ is said to be  \emph{graded} if it satisfies the condition $\calS=\oplus_{g\in G} \calS\cap\calL$. In particular, the derived ideal $\calL'=[\calL,\calL]$ is always graded, as $\calL'=\sum_{g,h\in G}[\calL_g,\calL_h]$ is the sum of graded subspaces.

Another grading $\calL=\oplus_{\gamma\in\Gamma}\calL_\gamma$ is said to be a \emph{coarsening} of the previous one in case for any $g\in G$ there is a $\gamma\in \Gamma$ such that $\calL_g\subseteq \calL_\gamma$. In other words, each $\calL_\gamma$ is a sum (necessarily direct) of homogeneous subspaces of the first grading. In this situation, the first grading is said to be a \emph{refinement} of the second one.

In what follows, we will work over an arbitrary ground field $\bF $ (even characteristic $2$ is allowed).

\medskip

\section{A counterexample of minimal dimension}

\noindent\textbf{The counterexample.}\quad Consider the four dimensional Lie algebra $\calL$ with a basis $\{a,u,v,w\}$ and multiplication given by:
\[
[a,u]=u,\quad [a,v]=w,\quad [a,w]=v,
\]
all the other brackets being $0$. Thus $\calL$ is the semidirect sum of the one dimensional subalgebra spanned by $a$ and the three dimensional abelian ideal spanned by $u$, $v$ and $w$. Then $\calL$ is graded as follows:
$\calL=\calL_\alpha\oplus\calL_\beta\oplus\calL_\gamma$, with
\[
\calL_\alpha=\bF a+\bF u,\quad \calL_\beta=\bF v,\quad \calL_\gamma=\bF w.
\]
If this were a semigroup grading, there would be a semigroup $\Gamma$ with three different elements $\alpha$, $\beta$ and $\gamma$ satisfying the following conditions (the binary operation in the semigroup will be denoted by juxtaposition, and note that it is not even necessary to assume that the semigroup is commutative):
\[
\begin{split}
\alpha^2&=\alpha,\quad\textrm{as $[\calL_\alpha,\calL_\alpha]=\bF u\subseteq \calL_\alpha$,}\\
\alpha\beta&=\gamma,\quad\textrm{as $[\calL_\alpha,\calL_\beta]=\bF w= \calL_\gamma$,}\\
\alpha\gamma&=\beta,\quad\textrm{as $[\calL_\alpha,\calL_\gamma]=\bF v= \calL_\beta$.}
\end{split}
\]
But then we would obtain:
\[
\gamma=\alpha\beta=\alpha^2\beta=\alpha(\alpha\beta)=\alpha\gamma=\beta,
\]
a contradiction. \qed

\medskip

The next result shows that the dimension of this counterexample is minimal.

\begin{theorem*} For any grading on a Lie algebra of dimension $\leq 3$ there exists a semigroup satisfying the property (P).
\end{theorem*}
\begin{proof}
The result is trivial for Lie algebras of dimension $1$ and very easy for dimension $2$, so it will be assumed that $\calL$ is a three dimensional graded Lie algebra, with a grading with two or three homogeneous subspaces: $\calL=\oplus_{\gamma\in\Gamma}\calL_\gamma$. Several possibilities may occur:
\begin{enumerate}
\item
If the dimension of the derived subalgebra $\calL'=[\calL,\calL]$ is $1$ then, since $\calL'$ is a graded subalgebra, there is an element $\alpha\in \Gamma$ such that $\calL'\subseteq \calL_\alpha$. Then $\Gamma$ becomes a (commutative) semigroup by means of $\gamma\delta=\alpha$ for any $\gamma,\delta\in\Gamma$, and the grading is obviously a grading over this semigroup.

\item
If the dimension of $\calL'$ is $2$ there are several subcases. Note that, as shown in \cite[p.~12]{J}, $\calL'$ is abelian and the center of $\calL$ is trivial.
\begin{enumerate}
\item There is an element $\alpha\in \Gamma$ such that $\calL'=\calL_\alpha$. But then $\calL=\calL_\alpha\oplus\calL_\beta$ for some other element $\beta\in\Gamma$, and $[\calL_\beta,\calL_\alpha]=\calL_\alpha$. This is then a grading over the integers, with $\calL_0=\calL_\beta$ and $\calL_1=\calL_\alpha$.
\item $\calL=\calL_\alpha\oplus\calL_\beta$ with $\dim\calL_\alpha=2$ and $\calL'=(\calL'\cap\calL_\alpha)\oplus\calL_\beta$. Here $\calL'=[\calL_\alpha,\calL_\alpha]+[\calL_\alpha,\calL_\beta]$ and $[\calL_\alpha,\calL_\beta]\ne 0$ as the center of $\calL$ is trivial. Hence either $[\calL_\alpha,\calL_\alpha]=\calL_\beta$ and $[\calL_\alpha,\calL_\beta]=\calL'\cap\calL_\alpha$, in which case this is a $\bZ/2\bZ$-grading with $\calL_\beta=\calL\subo$ and $\calL_\alpha=\calL\subuno$, or $[\calL_\alpha,\calL_\alpha]\subseteq \calL_\alpha$ and $[\calL_\alpha,\calL_\beta]=\calL_\beta$, and we get a $\bZ$-grading with $\calL_\alpha=\calL_0$ and $\calL_\beta=\calL_1$.
\item $\calL=\calL_\alpha\oplus\calL_\beta\oplus\calL_\gamma$, with all the homogeneous subspaces of dimension $1$ and $\calL'=\calL_\beta\oplus\calL_\gamma$. Since $\calL'$ is abelian, we get $\calL'=[\calL_\alpha,\calL_\beta]\oplus[\calL_\alpha,\calL_\gamma]$, so either $[\calL_\alpha,\calL_\beta]=\calL_\beta$ and $[\calL_\alpha,\calL_\gamma]=\calL_\gamma$, in which case we have a $\bZ$-grading with $\calL_\alpha=\calL_0$, $\calL_\beta=\calL_1$ and $\calL_\gamma=\calL_{-1}$, or $[\calL_\alpha,\calL_\beta]=\calL_\gamma$ and $[\calL_\alpha,\calL_\gamma]=\calL_\beta$, which is a $(\bZ/2\bZ)^2$-grading with $\calL_\alpha=\calL_{(\bar 1,\bar 1)}$, $\calL_\beta=\calL_{(\bar 1,\bar 0)}$ and $\calL_\gamma=\calL_{(\bar 0,\bar 1)}$.
\end{enumerate}
\item Finally, if $\calL'=\calL$ then $\calL$ is simple (see \cite[p.~12]{J}). Again there are two subcases:
\begin{enumerate}
\item Assume that $\calL=\calL_\alpha\oplus\calL_\beta$ with $\dim\calL_\alpha=1$, $\dim\calL_\beta=2$. Then $\calL=\calL'=[\calL_\alpha,\calL_\beta]+[\calL_\beta,\calL_\beta]$. Note that the first summand is not contained in $\calL_\alpha$, as this would force $\calL_\alpha$ to be an ideal. Hence we have $[\calL_\alpha,\calL_\beta]=\calL_\beta$, $[\calL_\beta,\calL_\beta]= \calL_\alpha$, and this is clearly a $\bZ/2\bZ$-grading.
\item Otherwise, $\calL=\calL_\alpha\oplus\calL_\beta\oplus\calL_\gamma$, for one dimensional subspaces $\calL_\alpha$, $\calL_\beta$ and $\calL_\gamma$. If there are two indices $\mu,\nu$ such that $[\calL_\mu,\calL_\nu]\subseteq \calL_\nu$, we may assume without loss of generality that $[\calL_\alpha,\calL_\beta]=\calL_\beta$. Hence $\calL_\alpha=\bF h$, $\calL_\beta=\bF x$ with $[h,x]=x$. Since the trace of the adjoint action by any element is $0$ (as $\calL=\calL'$), it follows that there exists an element $y\in\calL$ with $\calL_\gamma=\bF y$ and $[h,y]=-y$. Also, since $\calL=\calL'$, $[y,z]$ must belong to $\calL_\alpha$, and this gives a $\bZ$-grading with $\calL_\alpha=\calL_0$, $\calL_\beta=\calL_1$ and $\calL_\gamma=\calL_{-1}$. We are left with the case in which $[\calL_\alpha,\calL_\beta]=\calL_\gamma$, $[\calL_\beta,\calL_\gamma]=\calL_\alpha$, and $[\calL_\gamma,\calL_\alpha]=\calL_\beta$, and this is a $(\bZ/2\bZ)^2$-grading with trivial zero homogenous space. \qedhere
\end{enumerate}
\end{enumerate}
\end{proof}

\smallskip

\begin{remark*}
Consider the three dimensional Lie algebra $\calL$ with a basis $\{x,y,z\}$ and multiplication given by $[x,z]=[y,z]=z$, $[x,y]=0$. This is a direct sum of the two dimensional non abelian Lie algebra spanned by $y$ and $z$ and the one dimensional center spanned by $x-y$. The grading where the homogeneous spaces are the subspaces spanned by the basic elements is not a group grading (otherwise the indices of both $x$ and $y$ should correspond to the neutral element). This gives an example of minimal dimension of a non group grading on a Lie algebra.
\end{remark*}

\section{A counterexample on a semisimple Lie algebra}

\noindent\textbf{The counterexample.}\quad Consider now the three dimensional simple Lie algebras $\calJ=\espan{h,x,y}$ and $\calK=\espan{e_1,e_2,e_3}$, with multiplication given by:
\[
\begin{aligned}%
&[h,x]=x,&& [h,y]=-y,&& [x,y]=h,\\
&[e_1,e_2]=e_3,\ && [e_2,e_3]=e_1,\ && [e_3,e_1]=e_2.
\end{aligned}
\]
If the characteristic of the ground field is $\ne 2$, then $\calJ$ is isomorphic to $\frsl(2)$, and $\calK$ to the orthogonal Lie algebra $\frso(3)$. If, in addition, the ground field contains the square roots of $-1$, then both Lie algebras are isomorphic, and its direct sum is isomorphic to the orthogonal Lie algebra $\frso(4)$.

Let us take the Lie algebra $\calL=\calJ\oplus\calK$ with the grading given by:
\[
\calL_\alpha=\bF h+\bF e_1,\quad \calL_\beta=\bF x,\quad \calL_\gamma=\bF y,\quad
\calL_\delta=\bF e_2,\quad \calL_\mu=\bF e_3.
\]
It is straightforward to check that this is indeed a Lie grading. But if this were a semigroup grading, we would have the following condition in the semigroup:
\[
\alpha=\beta\gamma=(\alpha\beta)\gamma=\alpha(\beta\gamma)=\alpha^2,
\]
and, therefore,
\[
\mu=\alpha\delta=\alpha^2\delta=\alpha(\alpha\delta)=\alpha\mu=\delta,
\]
a contradiction.
\qed

\bigskip

Note that $\calJ$ is $\bZ$-graded with $\calJ_0=\bF h$, $\calJ_1=\bF x$, and $\calJ_{-1}=\bF y$. On the other hand, $\calK$ is $(\bZ/2\bZ)^2$-graded with $\calK_{(\bar 1,\bar 0)}=\bF e_1$, $\calK_{(\bar 0,\bar 1)}=\bF e_2$, and $\calK_{(\bar 1,\bar 1)}=\bF e_3$. Therefore, $\calL=\calJ\oplus\calK$ is $\bZ\times(\bZ/2\bZ)^2$-graded, with all the homogeneous spaces of dimension $1$. The grading in the counterexample above is a coarsening of this grading, obtained by joining together the homogeneous subspaces spanned by $h$ and $e_1$.

\end{document}